\documentclass[11pt, a4paper]{article}
\usepackage{graphicx,indentfirst}
\usepackage{amsmath,amssymb,mathrsfs}
\usepackage{amsthm}
\usepackage{verbatim}
\usepackage{appendix}
\usepackage{titletoc}
\usepackage{titlesec}
\usepackage{appendix}
\numberwithin{equation}{section}
\usepackage{amsfonts}
\usepackage[a4paper,left=2.5cm,right=2.5cm,bottom=3cm,top=3cm]{geometry}

\allowdisplaybreaks
 \theoremstyle{plain}\newtheorem{thm}{Theorem}[section]
\theoremstyle{plain}\newtheorem{lem}{Lemma}[section]
\theoremstyle{definition}
\theoremstyle{definition}
\theoremstyle{plain}
\theoremstyle{plain}
 \theoremstyle{remark}\newtheorem{rem}{Remark}[section]

\newcommand{\ds}{d}

\newcommand{\tr}{\mathrm{tr}}
\newcommand{\dv}{dv_{g_t}}
\newcommand{\td}{\tilde}
\newcommand{\vs}{v^{(6)}}
\begin{document}
\title{The Second Variational Formula For the Functional $\int v^{(6)}(g)dV_g$}
\author{Bin Guo \and Haizhong Li\thanks{Supported by grants of NSFC-10971110.}}
\date{}





\maketitle
\begin{abstract}
\noindent In this note, we compute the second variational formula for the functional
$\int_M v^{(6)}(g)dv_g$, which was introduced by Graham-Juhl [GJ] and the first
variational formula was obtained by Chang-Fang (\cite{CF08}). We also prove that Einstein
manifolds (with dimension $\ge 7$) with positive scalar curvature is a strict local
maximum within its conformal class, unless the manifold is isometric to round sphere with
the standard metric up to a multiple of constant. Note that when $(M,g)$ is locally
conformally flat, this functional reduces to the well-studied $\int_M \sigma_3(g)dv_g$.
Hence, our result generalize a previous result of Jeff Viaclovsky (\cite{V}) without the
locally conformally flat restraint.
\end{abstract}

\section{Introduction}
In the following, we let $(M^n,g)$ denote a compact, connected, smooth Riemannian
manifold without boundary. We denote the Ricci curvature and scalar curvature by $Ric$
and $R$, respectively. Recall that the Schouten tensor $P_{ij}$ is defined by
\begin{equation*}
P_{ij}=\frac{1}{n-2}\Big(R_{ij}-\frac{R}{2(n-1)}g_{ij}\Big),
\end{equation*}
and the Riemann curvature tensor can be written by
\begin{equation*}
Riem = W + P\odot g,
\end{equation*}
where $W$ is the Weyl curvature and $\odot$ is the Kulkarni-Nomizu product, which is
defined by
\begin{equation*}
(\alpha\odot\beta)_{ijkl}=\alpha_{ik}\beta_{jl}+\alpha_{jl}\beta_{ik}-\alpha_{il}\beta_{jk}-\alpha_{jk}\beta_{il},
\quad \forall\text{ symmetric 2-tensors }\alpha, \beta.
\end{equation*}
The $\sigma_k (g)$ curvature is defined to be the $k$-th elementary symmetric polynomial
of the eigenvalue of the Schouten tensor $P$. In [V], Viaclovsky started study of the
variational problems of the functional $\int_M \sigma_k(g)dv_g$, he proved that the first
variation of the functional $\int_M \sigma_k(g)dv_g (k=1,2)$ within a conformal class
subject to the constraint $Vol(M,g)=1$ is a metric satisfying $\sigma_k(g)\equiv
\text{const}$, and if $k\ge 3$ and the Riemannian manifold is locally conformally flat,
the same result follows. 
However, for $k\ge 3$ and the manifold is not locally conformally flat,
$\sigma_k(g)\equiv \text{const}$ is not Euler-Lagrange equation of the functional $\int_M
\sigma_k dv$ within a conformal class subject to the constraint $Vol(M,g)=1$.

The renormalized volume coefficients of $g$, denoted here by $v^{(2k)}(g)$, arose in the
late 90s in the physics literature. They are defined in terms of the expansion of the
ambient or Poincare metric associated to $g$. If the Riemannian manifold is locally
conformally flat, these quantities coincide with the $\sigma_k(g)$ up to a constant. More
precisely, it is known that (see [GJ], [CF], or [GHL])
\begin{equation*}
\begin{split}
 v^{(2)}(g)=-\frac{1}{2}\sigma_1(g),\qquad
 v^{(4)}(g)=\frac{1}{4}\sigma_2(g),\\
  v^{(6)}(g)=-\frac{1}{8}\big[\sigma_3(g)+\frac{1}{3(n-4)}(P_g)^{ij}(B_g)_{ij}\big],
\end{split}
\end{equation*}
 where
\begin{equation}\label{eqn:ba}
 (B_g)_{ij}:=\frac{1}{n-3}\nabla^k\nabla^lW_{likj}+\frac{1}{n-2}R^{kl}W_{likj}
\end{equation}
 is the {\it Bach} tensor of the metric.
 Just as
$\int_M \sigma_k(g^{-1}\circ A_g)\,dv_g$ is conformally invariant when $2k=n$ and $(M,g)$
is locally conformally flat, Graham  showed in \cite{G00} that $\int_M v^{(2k)}(g)
\,dv_g$ is also conformally invariant on a general manifold when $2k=n$. Chang and Fang
showed in \cite{CF08} that, for $n\neq 2k$, the Euler-Lagrange equations for the
functional $\int_M v^{(2k)}(g) \,dv_g$ under conformal variations subject to the
constraint $Vol_g(M)=1$ satisfies $v^{(2k)}(g)=$ const., which is a generalized
characterization for the curvatures $\sigma_k(g^{-1}\circ A_g)$ when $(M,g)$ is locally
conformally flat, as given by Viaclovsky \cite{V}.

We note that Graham \cite{G00} also gives an explicit expression of $v^{(8)}(g)$, but the
explicit expression of $v^{(2k)}(g)$ for general $k$ is not known because they are
algebraically complicated (see page 1958 of \cite{G00}). Thus the study of the
$v^{(2k)}(g)$ curvatures involves significant challenges not shared by that of
$\sigma_k(g)$:  firstly, for $k\geq 3$, $v^{(2k)}(g)$ depends on derivatives of curvature
of $g$--- in fact, for $k\geq 3$, $v^{(2k)}(g)$ depends on derivatives of curvatures of
order up to $2k-4$; secondly, the $v^{(2k)}(g)$ are defined via an indirect highly
nonlinear inductive algorithm (see \cite{G00}).
 We aim to study the stability of the critical metric of the functional
\begin{equation*}
\mathcal{F}_3[g]=\frac{\int _M v^{(6)}(g)dv_g}{\big(\int_M dv_g\big)^{(n-6)/n}},
\end{equation*}
within a conformal class. First we recall the theorem of Chang-Fang [CF] (also see Graham
[G]).
\begin{thm}\label{thm:thm2}(\cite{CF08})
Let $(M^n,g)$ be an $n$-dimensional $(n\geq 7)$ compact Riemannian manifold, then the
functional $\mathcal{F}_3(g)$ is variational within the conformal class, i.e. the
critical metric in $[g]$ satisfies the equation
\begin{equation}\label{eqn:first}
v^{(6)}\equiv \text{const}.
\end{equation}
If $n=6$, $\mathcal{F}_3[g]$ is a constant in the conformal class $[g]$.
\end{thm}
In this note, we compute the second variational formula of $\mathcal{F}_3[g]$ within its
conformal class $[g]$. Our results are
\begin{thm}\label{thm:thm1}
Let $(M^n,g)$ be an $n$-dimensional $(n\geq 7)$  compact Riemannian manifold with
$\vs(g)=$const, then the second variational formula of the functional
$\mathcal{F}_3{[g]}$ within its conformal class at $g$ is
\begin{equation*}
\frac{d^2}{dt^2}\Big|_{t=0} \mathcal{F}_3[g_t]= (n-6)V^{-(n-6)/n}\Bigg\{\int
\Big[-6\vs(g)\bar{\phi}^2+\Big(
\frac{B_{ij}\bar{\phi}_{ij}}{24(n-4)}+\frac{1}{8}T_{2ij}\bar{\phi}_{ij}-\frac{1}{12}P_{ij}C_{ijk}\bar{\phi}_k\Big)\bar{\phi}
\Big]dv\Bigg\}
\end{equation*}
where $g_t=e^{2u_t}g$, $\frac{\partial}{\partial t}\big|_{t=0}u_t=\phi$, and
$\bar\phi=\phi-\frac{\int_M\phi dv_g}{\int_M dv_g}$, $T_{2ij}$ and $C_{ijk}$ are defined
in section 2.
\end{thm}

\begin{thm}\label{thm:main}
Let $(M^n,g)$ be an $n$-dimensional $(n\geq 7)$  compact Einstein manifold with positive
scalar curvature. Then it is a strict local maximum within its conformal class $[g]$,
unless $(M^n,g)$ is isometric to $S^n$ with the standard metric up to a multiple of
constant.
\end{thm}

\begin{rem} When $(M^n,g)$ is a locally conformally flat, $v^{(6)}(g)=-\frac{1}{8}\sigma_3(g)$, for
the functional $\int_M \sigma_3(g)dv_g$, J. Viaclovsky (\cite{V}) proved that a positive
constant sectional curvature metric is a strict local minimum, unless the manifold is
isometric to $S^n$ with the standard metric. Our result coincides with his at the locally
conformally flat Einstein metrics, however, ours does not need the locally conformally
flat assumption.
\end{rem}

\section{Preliminaries}
Let $(M^n,g)$ be an $n$-dimensional compact Riemannian manifold. Throughout this note, we
make the convention that repeated index means summation over $1$ to $n$. First we recall
the transformation law of various curvatures under conformal change of metrics. Let
$\td{g}=e^{2u}g$, $u\in C^\infty (M)$, then the Riemannian curvature tensors satisfy
\begin{equation*}
Riem(\td{g})=e^{2u}(Riem (g) - \alpha\odot g),
\end{equation*}
where $\alpha_{ij}=u_{ij}-u_iu_j+\frac{|\nabla u|^2}{2}g_{ij}$ (note that $u_{ij}$ means
the covariant derivative with respect to the fixed metric $g$). By contracting, we see
that the Ricci curvature and scalar curvature satisfy
\begin{equation}\label{eqn:trans}
R_{ij}(\td{g})= R_{ij}-(n-2)\alpha_{ij}-(\sum_k\alpha_{kk})g_{ij}, \quad R(\td{g})=
e^{-2u}R - 2(n-1)e^{-2u}\sum_k\alpha_{kk}.
\end{equation}
From \eqref{eqn:trans} and the definition of Schouten tensor, we see that
\begin{equation}\label{eqn:schouten}
\td{P}_{ij}=P_{ij}-\alpha_{ij},
\end{equation}
where we denote $P(\td{g})$ by $\td{P}$ for notations convenience.

\begin{lem}\label{lem:lem1}
We have the following formulae (see e.g. [GHL])

(1) $\nabla^i W_{ijkl}=-(n-3)C_{jkl}, \; C_{ijk}\text{ is the {\it{Cotton}} tensor
defined by } P_{ij,k}-P_{ik,j}$;

(2) \label{eqn:bach} $\nabla^j B_{ij}=(n-4)\sum_{k,l}P_{kl}C_{kli}$;

(3) $B_{ij}=B_{ji}$, where $B_{ij}$ is defined by \eqref{eqn:ba}.

\end{lem}
The proof of Lemma \ref{lem:lem1} is a direct calculation and one can find it in
\cite{GHL}.

Let $V$ be a vector space, $A:V\to V$ a linear map. Define the Newton transformation
$T_k(A)$ $:V\to V$ by:
\begin{equation*}
T_{k}(A):=\sigma_k(A)I-\sigma_{k-1}(A)A+\cdots+(-1)^k
A^k=\sum_{i=0}^k\sigma_{k-i}(A)(-1)^iA^i,
\end{equation*}
where $I$ is the identity map and $\sigma_k(A)$ is the $k$-th elementary symmetric
polynomial of the eigenvalues of $A$. Under an orthnormal basis of $V$, $T_{k}$ can be
written as follows:
\begin{equation*}
{T_{k}}_{ij}=\frac{1}{k!}\delta^{j_1\ldots j_k j}_{i_1\ldots i_k i}A_{i_1j_1}\ldots
A_{i_kj_k},
\end{equation*}
where $\delta^{j_1\ldots j_k j}_{i_1\ldots i_k i}$ is the generalized Kronecker notation.
We  recall some well-known results in this respect, which we will need in our later
arguments.
\begin{lem}\label{lem:newton} The Newton transformations
$T_{k}$ satisfy ([R], [GHL])

(1) Newton's formula: $(k+1)\sigma_{k+1}(A)=\tr(T_kA)$;

(2) $\frac{d}{dt}\sigma_k(A_t)=\tr(T_{k-1}\frac{d}{dt} A_t)$, for any family of
transformations $A_t: V\to V$.

(3) $\tr(T_{k})=(n-k)\sigma_k(A)$.

\end{lem}

In the following we denote $T_k(g^{-1}\circ P)$ simply by $T_k$. We have the following
formula, which is a direct calculation (see [GHL])
\begin{equation}\label{eqn:tij}
\sum_{i}{T_2}_{ij,i}=-\sum_{k,l}P_{kl}C_{klj}.
\end{equation}

\section{The first variational formula and proof of Theorem 1.1}
In this section, we will compute the Euler-Lagrange equation for the functional
$\mathcal{F}_3(g)$ within the conformal class. For convenience we denote the numerator of
$\mathcal{F}_3(g)$ by

$$
F(g)=\int_M v^{(6)}(g)dv_g.
$$
 Under the conformal change of metrics
$g_t=e^{2u(t)}g$, by use of \eqref{eqn:schouten}, we see that in local coordinates (see
[CF])
\begin{equation*}
\tilde{P}_{i}^{\;j}=e^{-2u(t)}(P_{i}^j-\alpha_{i}^{\;j})
\end{equation*}
\begin{equation*}
\tilde{B}_{i}^{\;j}=e^{-4u(t)}\Big(B_{i}^{\;j}+(n-4)u^k\big(g^{jl}C_{ilk}+g^{jl}C_{lik}\big)+(n-4)u^k
u^lg^{pj}W_{ikpl}\Big),
\end{equation*}
where we write $\alpha_{ij}=u_{ij}(t)-u_i(t)u_j(t)+\frac{1}{2}|\nabla_g u(t)|^2g_{ij}$
and we make the convention that $\tilde{P}_{ij} = P_{ij}(g_t)$,
$\tilde{B}_{ij}=B_{ij}(g_t)$, etc, for notations convenience.

For notions convenience we denote $\frac{d}{d t}$ by $\delta$. Denote
$\frac{\partial}{\partial t}\Big|_{t=0}u=\phi$, and $\frac{\partial^2}{\partial
t^2}\Big|_{t=0}u=\psi$. With the above preparations, we have
\begin{equation*}
\delta \tilde{P}_{i}^{\;j}=-2(\delta u)\tilde{P}_{i}^{\;j}-e^{-2u}\delta\alpha_{i}^{\;j}
\end{equation*}
\begin{equation*}
\begin{split} \delta \tilde{B}_{i}^{\;j}=-4(\delta u)\tilde{B}_{i}^{\;j}+(n-4)e^{-4u}\Big(&(\delta
u)^k\big(g^{jl}C_{ilk}+g^{jl}C_{lik}\big)\\
        & +(\delta u)^k u^lg^{pj}W_{ikpl}+u^k(\delta u)^lg^{pj}W_{ikpl}\Big)
\end{split}\end{equation*}
 \begin{equation*}
\delta (dv_{g_t})= n(\delta u) dv_{g_t}.
\end{equation*}
Now we derive the first variation formula for $\mathcal{F}_3$. First we have

\begin{equation}\label{eqn:new1st}
-8\delta F=\int
\delta(\sigma_3(g_t)dv_{g_t})+\frac{1}{3(n-4)}\delta(\tilde{P}_{i}^{\;j}\tilde{B}_{j}^{\;i}dv_{g_t}).
\end{equation}
By use of Lemma \ref{lem:newton}, we have
\begin{align*}
\delta(\sigma_3(g_t))=\tilde{T}_{2j}^{\;i}\delta\tilde{P}_{i}^{\;j}&=\tilde{T}_{2j}^{\;i}
                           (-2(\delta u)\tilde{P}_{i}^{\;j}-e^{-2u}\delta\alpha_{i}^{\;j})\\
                  &=-6(\delta u)\sigma_3(g_t)-e^{-2u}{\tilde{T}_{2j}}^{\;i}\delta\alpha_{i}^{\;j}.
\end{align*}
Hence
\begin{equation}\label{eqn:11}\begin{split}
\int \delta(\sigma_3(g_t)dv_{g_t})&=\int \Big[-6(\delta u)\sigma_3(g_t)
                      -e^{-2u}\tilde{T}_{2i}^{\;j}\delta\alpha_{j}^{\;i}+n(\delta u)\sigma_3(g_t)\Big]dv_{g_t}\\
                      &=\int \Big[(n-6)(\delta u)\sigma_3(g_t)-e^{-2u}\tilde{T}_{2i}^{\;j}\delta\alpha_{j}^{\;i}\Big]dv_{g_t}.
\end{split}\end{equation}
On the other hand, the second term of \eqref{eqn:new1st} is
\begin{align}
\label{eqn:22}
&\int \frac{1}{3(n-4)}\delta(\tilde{P}_{i}^{\;j}\tilde{B}_{j}^{\;i}dv_{g_t})\\
=& \int \frac{1}{3(n-4)}\bigg[\delta\tilde{P}_{i}^{\;j}\tilde{B}_{j}^{\;i}
    +\tilde{P}_{i}^{\;j}\delta\tilde{B}_{j}^{\;i}+n(\delta u)\tilde{P}_{i}^{\;j}\tilde{B}_{j}^{\;i}\bigg]\dv\notag\\
     = &\int \frac{1}{3(n-4)}\bigg[\tilde{B}_{j}^{\;i}
\Big((-2\delta u)\tilde{P}_{i}^{\;j}-e^{-2u}(\delta\alpha_{i}^{\;j})\Big)+
     \tilde{P}_{i}^{\;j}\Big(-4(\delta u)\tilde{B}_{j}^{\;i} + (n-4)e^{-4u}\big((\delta u)^k(g^{il}C_{ljk}+C_{jlk})\notag \\
                                                &\qquad\quad+2(\delta u)^k\ u^lg^{ip}W_{pkjl}\big)\Big)
                                                +n(\delta u)\tilde{P}_{i}^{\;j}\tilde{B}_{j}^{\; i}
                                             \bigg]\dv\notag\\
     &=\int \frac{1}{3(n-4)}\bigg[(n-6)(\delta u)\tilde{P}_{i}^{\;j}\tilde{B}_{j}^{\;i}
                  -e^{-2u}\tilde{B}_{j}^{\;i}
                  (\delta\alpha_{i}^{\;j})\notag \\ &\qquad
      +2(n-4)e^{-4u}\tilde{P}_{i}^{\;j}\Big((\delta u)^kg^{il}C_{ljk}+(\delta u)^k
      u^lg^{ip}W_{pkjl}\Big)\bigg]\dv\notag
\end{align}

From calculations in \eqref{eqn:11} and \eqref{eqn:22}, we have the following formula,
which will be used in section 4
\begin{align}\label{eqn:vv}
&\delta \big(v^{(6)}(g_t)\dv\big)\\
=&-\frac{1}{8}\Big[(n-6)(\delta
u)\sigma_3(g_t)-e^{-2u}{\tilde{T}}_{2i}^{\;j}(\delta\alpha)_{j}^{\;i}\Big]\dv
  -\frac{1}{24(n-4)}\Big[(n-6)(\delta u)\tilde{P}_{i}^{\;j}\tilde{B}_{j}^{\;i}\notag\\ &\;\quad
      -e^{-2u}\tilde{B}_{i}^{\;j}(\delta\alpha)_{j}^{\;i}
                       +2(n-4)e^{-4u}\tilde{P}_{i}^{\;j}\Big((\delta u)^kg^{il}C_{ljk} +
                         (\delta u)^k
                         u^lg^{ip}W_{pkjl}\Big)\Big]\dv.\notag
\end{align}

Thus we have
\begin{align}\label{eqn:1st}
\delta F=&\int \delta (v^{(6)}(g_t)\dv)\\
   =&-\frac{1}{8}\int\Bigg\{(n-6)(\delta u)\sigma_3(g_t)-e^{-2u}\tilde{T}_{2i}^{\;j}\delta\alpha_{j}^{\;i}
 + \frac{1}{3(n-4)}\Big[(n-6)(\delta u)\tilde{P}_{i}^{\;j}\tilde{B}_{j}^{\;i}\notag\\ &\;
 -e^{-2u}\tilde{B}_{j}^{\;i}(\delta\alpha)_{i}^{\;j}
               +2(n-4)e^{-4u}\tilde{P}_{i}^{\;j}\big((\delta u)^kg^{il}C_{ljk}+(\delta u)^k u^lg^{ip}W_{pkjl}\big)\Big]\Bigg\}\dv\notag\\
 =&\int \Bigg\{(n-6)(\delta u)v^{(6)}(g_t)-\frac{1}{8}\bigg[ -e^{-2u}\tilde{T}_{2i}^{\;j}\delta\alpha_{j}^{\;i} -
                \frac{1}{3(n-4)}e^{-2u}\tilde{B}_{j}^{\;i}
 \delta\alpha_{i}^{\;j}\notag\\
 & +\frac{2}{3}e^{-4u}\tilde{P}_{i}^{\;j}\Big((\delta u)^k g^{il}C_{ljk}+
        (\delta u)^k u^lg^{ip} W_{pkjl}\Big)\bigg]\Bigg\}dv_{g_t}.\notag
\end{align}

\noindent {\bf Proof of Theorem 1.1} Noting that $u(0)=0$, we conclude the first
variational formula of $\mathcal{F}_3[g_t]$ within the conformal class $[g]$ is (see [CF]
or [G])
\begin{align}\label{eqn:main1}
\frac{d}{dt}\Big|_{t=0}\mathcal{F}_3(g_t)&=\frac{d}{dt}\Big|_{t=0}F\cdot
V^{-\frac{n-6}{n}}-
\frac{n-6}{n}V^{-\frac{n-6}{n}}\Big(\int v^{(6)}\ds v\Big)\int n \phi \ds v\\
=&V^{-\frac{n-6}{n}}\bigg\{(n-6)\int \phi v^{(6)}(g)+\frac{1}{8}\int
T_{2ij}\phi_{ij}+\frac{1}{8}\int \frac{B_{ij}\phi_{ij}}{3(n-4)}
  -\frac{1}{12}
  \int P_{ij}\phi_kC_{ijk}\notag \\ &\qquad \; \qquad-(n-6)V^{-1}\Big(\int v^{(6)}(g)\Big)\int \phi \ds v\bigg\}\notag\\
=&(n-6)V^{-\frac{n-6}{n}}\bigg\{\int \phi \Big(v^{(6)}-V^{-1}\int v^{(6)}\Big)\ds
v\bigg\}.\notag
\end{align}
where we have used \eqref{eqn:tij} and (2) of Lemma \ref{lem:lem1} and the integration by
parts. Here $V=\int \ds v_g$. Hence, we see that the Euler-Langrange equation of the
functional $\mathcal{F}_{3}(g)$ within the conformal class $[g]$ is
\begin{equation*}
\vs(g)=V^{-1}\int \vs(g)dv_g\equiv \text{const},
\end{equation*}
and we get Theorem \ref{thm:thm2}.\qed

\section{The Second Variational Formula and proofs of Theorem 1.2-1.3}
In this section, we will calculate the second variational formula for the functional
$\mathcal{F}_3$ within the conformal class $[g]$. The computation is direct and routine.
For convenience, we separate each term in the first variational equation \eqref{eqn:1st}
and compute them respectively.

For derivative of the first term in \eqref{eqn:1st}, by use of \eqref{eqn:vv}, we have
\begin{align}\label{eqn:m1}
&(n-6)\frac{d}{dt}\Big|_{t=0}\int (\delta u) v^{(6)}(g_t)\dv\\
=& (n-6)\int\Bigg\{ \psi v^{(6)}(g)-\frac{1}{8}(\delta u) \Big[(n-6)(\delta
u)\sigma_3(g_t)-
                          e^{-2u} \tilde{T}_{2i}^{\;j}\delta\alpha_{j}^{\;i}\Big]\dv\Big|_{t=0}\notag\\
  &\quad\qquad\qquad-\frac{(\delta u)}{24(n-4)}\Big[(n-6)(\delta u)\tilde{P}_{i}^{\;j}\tilde{B}_{j}^{\;i}
      - e^{-2u}\tilde{B}_{j}^{\;i}\delta\alpha_{i}^{\;j}\notag
        \\ &\quad\qquad\qquad +2(n-4)e^{-4u}\tilde{P}_{i}^{\;j}\bigg((\delta u)^kg^{il}C_{ljk}
                   +(\delta u)^k u^lg^{ip}W_{pkjl}\bigg)\Big]\dv\Big|_{t=0}\Bigg\}\notag\\
  =& (n-6)\int \bigg\{\psi v^{(6)}(g)+(n-6)\phi^2v^{(6)}(g)-\frac{1}{8} \Big[-\phi T_{2ij}\phi_{ij}
          \notag \\ &\;\qquad -\frac{\phi B_{ij}\phi_{ij}}{3(n-4)}
  +\frac{2}{3}
                \phi\phi_k P_{ij} C_{ijk}\Big]\bigg\}dv.\notag
\end{align}
For derivative of the second term in \eqref{eqn:1st}, we need the following formula of
the variation of the Newton transformation:

\begin{align}
\frac{d}{dt}\Big|_{t=0}(\td{T}_{2i}^{\;j})&=\frac{d}{dt}\Big|_{t=0}\Big(\frac{1}{2!}\delta^{j_1j_2j}_{i_1i_2i}
\td{P}_{j_1}^{\;i_1}\td{P}_{j_2}^{\;i_2}\Big)
             =\delta^{j_1j_2j}_{i_1i_2i}P_{j_1}^{i_1}\frac{d}{dt}\Big|_{t=0}\td{P}_{j_2}^{i_2}\\
             &=\delta^{j_1j_2j}_{i_1i_2i}P_{i_1j_1}(-2\phi P_{i_2j_2}-\phi_{i_2j_2})=
             -4\phi {T_{2}}_{ij}-\delta^{j_1j_2j}_{i_1i_2i}P_{i_1j_1}\phi_{i_2j_2}.\notag
\end{align}
Therefore, the variation of the second term of \eqref{eqn:1st} is given by 
\begin{align}\label{eqn:m2}
\frac{1}{8}&\frac{d}{dt}\Big|_{t=0}\int
e^{-2u}\td{T}_{2i}^{\;j}\delta\alpha_{j}^{i}\dv=\frac{1}{8}\int \Bigg\{
      -2\phi
      T_{2ij}\phi_{ij}+\frac{d}{dt}\Big|_{t=0}\td{T}_{2i}^{\;j}\phi_{j}^{\;i}\\
&\quad\qquad +T_{2ij}\frac{d^2}{dt^2}\Big|_{t=0}\alpha_{ij}
      +  n\phi T_{2ij}\phi_{ij}\Bigg\}\ds v\notag\\
=&\frac{1}{8}\int \bigg\{(n-2)\phi
T_{2ij}\phi_{ij}+T_{2ij}\frac{d^2}{dt^2}\Big|_{t=0}\alpha_{ij}
   +\phi_{ij}\Big(-4\phi {T_{2}}_{ij}
       -\delta^{j_1j_2j}_{i_1i_2i}P_{i_1j_1}\phi_{i_2j_2}\Big)\bigg\}dv\notag\\
=&\int\bigg\{ \frac{(n-6)}{8}\phi
T_{2ij}\phi_{ij}+\frac{1}{8}T_{2ij}\frac{d^2}{dt^2}\Big|_{t=0}\alpha_{ij}
   -\frac{1}{8}\delta^{j_1j_2j}_{i_1i_2i}P_{i_1j_1}\phi_{i_2j_2}\phi_{ij}\bigg\}.\notag
\end{align}
The variation of the third term of \eqref{eqn:1st} is 
\begin{align}\label{eqn:m3}
&\frac{1}{24(n-4)}\frac{d}{dt}\Big|_{t=0}\int \Big(e^{-2u}\td{B}_{i}^{\;j}\delta\alpha_{j}^{\;i}\Big)\dv\\
=&\frac{1}{24(n-4)}\int\Bigg\{-2\phi
B_{ij}\phi_{ij}+\frac{d}{dt}\Big|_{t=0}(\td{B}_{i}^{\;j})\phi_{ij}
+B_{ij}\frac{d}{dt}\Big|_{t=0}\delta\alpha_{ij}+n\phi B_{ij}\phi_{ij}\Bigg\}\ds v\notag\\
=&\frac{1}{24(n-4)}\int\Bigg\{(n-2)\phi
B_{ij}\phi_{ij}+B_{ij}\frac{d^2}{dt^2}\Big|_{t=0}\alpha_{ij}
+\phi_{ij}\Big(-4\phi B_{ij}+2(n-4)\phi_kC_{ijk}\Big)\Bigg\}\ds v\notag\\
=&\int \Bigg\{\frac{n-6}{24(n-4)}\phi
B_{ij}\phi_{ij}+\frac{B_{ij}}{24(n-4)}\frac{d^2}{dt^2}\Big|_{t=0}\alpha_{ij}
+\frac{1}{12}\phi_k\phi_{ij}C_{ijk}\Bigg\}dv.\notag
\end{align}

The variation of the fourth term of \eqref{eqn:1st} is
\begin{align}\label{eqn:m4}
&-\frac{1}{12}\frac{d}{dt}\Big|_{t=0}\int \bigg[e^{-4u}\td{P}_{j}^{\;i}\Big((\delta
u)^kg^{jl}C_{ilk}+
     (\delta u)^k u^lg^{jp}W_{ikpl}\Big)\bigg]\dv\\
&=-\frac{1}{12}\int \Bigg\{-4\phi P_{ij}\phi_kC_{ijk} +
   \frac{d}{dt}\Big|_{t=0}\td{P}_{i}^{j}\phi^kg^{il}C_{ljk}
 +P_{ij}\psi_k C_{ijk}+P_{ij}\phi_k \phi_l W_{ikjl}+n\phi P_{ij}\phi_kC_{ijk}\Bigg\}dv\notag\\
&=-\frac{1}{12}\int \Bigg\{(n-4)\phi P_{ij}\phi_kC_{ijk}+P_{ij}\phi_k \phi_lW_{ikjl}
+\phi_kC_{ijk}\Big(-2\phi P_{ij}-\phi_{ij}\Big)+P_{ij}\psi_kC_{ijk}\Bigg\}dv\notag\\
&=\int \Bigg\{-\frac{(n-6)}{12}\phi
P_{ij}\phi_kC_{ijk}-\frac{1}{12}P_{ij}\phi_k\phi_lW_{ikjl}
+\frac{1}{12}\phi_k\phi_{ij}C_{ijk}-\frac{1}{12}P_{ij}\psi_kC_{ijk}\Bigg\}dv.\notag
\end{align}
Combining \eqref{eqn:m1}, \eqref{eqn:m2}, \eqref{eqn:m3} and \eqref{eqn:m4}, we have
\begin{align}\label{eqn:m2nd}
& \frac{d^2}{dt^2}\Big|_{t=0} F(g_t) \\ 
&=\int \Bigg\{(n-6)\psi v^{(6)}(g)+(n-6)^2\phi^2 v^{(6)}(g) +\frac{n-6}{4}\phi
T_{2ij}\phi_{ij}
-\frac{1}{8}\Big[ 
- \frac{2(n-6)\phi B_{ij}\phi_{ij}}{3(n-4)}\notag\\
 & \qquad+ \frac{4(n-6)}{3}\phi\phi_k P_{ij}C_{ijk}- T_{2ij}\frac{d^2}{dt^2}\Big|_{t=0}\alpha_{ij}-
 \frac{B_{ij}}{3(n-4)}\frac{d^2}{dt^2}\Big|_{t=0}\alpha_{ij}\notag\\
 &\qquad +\delta^{j_1j_2j}_{i_1i_2i}P_{i_1j_1}\phi_{i_2j_2}\phi_{ij}-\frac{4}{3}\phi_k\phi_{ij}C_{ijk}
 +\frac{2}{3}P_{ij}\phi_k\phi_l W_{ikjl}+\frac{2}{3}P_{ij}\psi_kC_{ijk}\Big]\Bigg\}dv\notag.
\end{align}
Since $\frac{d^2}{dt^2}\Big|_{t=0}\alpha_{ij}=\psi_{ij}-2\phi_i\phi_j+|\nabla
\phi|^2g_{ij}$,
 by use of divergence theorem we obtain
\begin{align}\label{eqn:alp}
&-\int T_{2ij}\frac{d^2}{dt^2}\Big|_{t=0}\alpha_{ij}-\int \frac{B_{ij}}{3(n-4)}\frac{d^2}{dt^2}\Big|_{t=0}\alpha_{ij}\\
&=-\int T_{2ij}(\psi_{ij}-2\phi_i\phi_j+|\nabla \phi|^2g_{ij})
    -\int \frac{B_{ij}(\psi_{ij}-2\phi_i\phi_j+|\nabla \phi|^2g_{ij})}{3(n-4)}\notag\\
&=\int
-\frac{2}{3}P_{kl}C_{kli}\psi_i+2T_{2ij}\phi_i\phi_j+\frac{2}{3(n-4)}B_{ij}\phi_i\phi_j
                        -|\nabla \phi|^2T_{2kk}\notag\\
&=\int -\frac{2}{3}P_{kl}C_{kli}\psi_i-|\nabla \phi|^2T_{2kk}-2T_{2ij,j}\phi_i\phi
   -2T_{2ij}\phi_{ij}\phi-\frac{2}{3(n-4)}B_{ij,j}\phi_i\phi-\frac{2}{3(n-4)}B_{ij}\phi_{ij}\phi\notag\\
&=\int -\frac{2}{3}P_{kl}C_{kli}\psi_i-|\nabla \phi|^2T_{2kk}+\frac{4}{3}\phi \phi_i
P_{kl}C_{kli}-2T_{2ij}\phi_{ij}\phi
   -\frac{2}{3(n-4)}B_{ij}\phi_{ij}\phi,\notag
\end{align}
where we have used the following identity in the second equality
$$
\int T_{2ij}\psi_{ij}dv+\frac{1}{3(n-4)}\int B_{ij}\psi_{ij}dv=\int \frac{2}{3}P_{kl}C_{kli}\psi_i,
$$
which can be checked by use of (2.3), (2) of Lemma 2.1 and integration by parts.

Substituting \eqref{eqn:alp} into \eqref{eqn:m2nd} and making some cancelations, we
conclude that
\begin{align}\label{eqn:2nd}
\frac{d^2}{dt^2}\Big|_{t=0} F(g_t)=&\int \Bigg\{(n-6) \psi
v^{(6)}+(n-6)^2\phi^2v^{(6)}(g)-\frac{1}{8}\Big[-2(n-5)\phi T_{2ij}\phi_{ij}
\\ &-\frac{2(n-5)}{3(n-4)}\phi B_{ij}\phi_{ij}+\frac{4(n-5)}{3}\phi P_{ij}\phi_kC_{ijk}\notag\\
  & -\frac{4}{3}\phi_k\phi_{ij}C_{ijk}+\delta^{mnj}_{kli}P_{km}\phi_{ln}\phi_{ij}
  +\frac{2}{3}P_{ij}\phi_k\phi_l W_{ikjl}
                     -(n-2)|\nabla \phi|^2\sigma_2(g)\Big]\Bigg\}dv\notag,
\end{align}
where we have used the identity that $T_{2kk}=(n-2)\sigma_2(g)$ (see Lemma 2.2). It
remains to study the last four terms on the right hand side of \eqref{eqn:2nd}.  By
definition,
\begin{align*}
\delta^{mnj}_{kli}&=\det\begin{pmatrix}\delta_{km}&\delta_{kn}&\delta_{kj}\\
                                       \delta_{lm}&\delta_{ln}&\delta_{lj}\\
                                       \delta_{im}&\delta_{in}&\delta_{ij}\end{pmatrix}\\
         =&\delta_{km}\delta_{ln}\delta_{ij}-\delta_{km}\delta_{lj}\delta_{in}-\delta_{lm}\delta_{kn}\delta_{in}+
           \delta_{lm}\delta_{in}\delta_{kj}+\delta_{im}\delta_{kn}\delta_{lj}-\delta_{im}\delta_{ln}\delta_{kj}.
\end{align*}
We compute by use of divergence theorem
\begin{equation}\label{eqn:mid}\begin{split}
&\int \delta^{mnj}_{kli}P_{km}\phi_{ln}\phi_{ij}=\int \phi\Big(\delta^{mnj}_{kli}(P_{mk}\phi_{nl})_{,ij}\Big)\\
=&\int
\phi\delta^{mnj}_{kli}\Big(P_{km,ij}\phi_{nl}+2P_{km,i}\phi_{nl,j}+P_{km}\phi_{nl,ij}\Big).
\end{split}\end{equation}
Now we compute integrands of the right hand side of \eqref{eqn:mid} respectively. The
first term is
\begin{align}\label{eqn:mid1}
&\phi\delta^{mnj}_{kli}P_{km,ij}\phi_{nl}\\ =&\phi
P_{km,ij}\phi_{nl}(\delta_{km}\delta_{ln}\delta_{ij}-\delta_{km}\delta_{lj}\delta_{in}-\delta_{lm}\delta_{kn}\delta_{in}+
           \delta_{lm}\delta_{in}\delta_{kj}+\delta_{im}\delta_{kn}\delta_{lj}-\delta_{im}\delta_{ln}\delta_{kj})\notag \\
           =&\phi P_{kk,ii}\phi_{nn}-\phi P_{kk,ij}\phi_{ij}-\phi P_{kl,ii}\phi_{kl}+\phi P_{kl,ik}\phi_{il}+\phi P_{ki,il}\phi_{kl}-\phi P_{ki,ik}\phi_{ll}\notag\\
           =&\phi \phi_{nn}C_{iik,k}+\phi\phi_{kl}C_{lik,i}+\phi\phi_{kl}C_{iki,l}
           =\phi \phi_{kl}C_{lik,i}.\notag
\end{align}
The second one is
\begin{align}\label{eqn:mid2}
&2\phi\delta^{mnj}_{kli}P_{km,i}\phi_{nlj}\\
=&2\phi P_{kk,i}\phi_{nni}-2\phi P_{kk,i}\phi_{ijj}-2\phi P_{km,i}\phi_{kmi}+2\phi P_{km,i}\phi_{imk}+2\phi P_{km,m}\phi_{kll}-2\phi P_{km,m}\phi_{llk}\notag\\
=&2\phi P_{km,i}(\phi_{mik}-\phi_{mki})=2\phi P_{km,i}\phi_{j}R_{jmik}.\notag
\end{align}
The third one is
\begin{align}\label{eqn:mid3}
&\phi\delta^{mnj}_{kli}P_{km}\phi_{nl,ij}\\
=&\phi \phi_{kk,ii}P_{nn}-\phi\phi_{iikl}P_{kl}-\phi\phi_{klii}P_{kl}+\phi\phi_{ilki}P_{kl}+\phi\phi_{kiil}P_{kl}-\phi\phi_{kiik}P_{ll}\notag \\
=&\phi P_{nn}(\phi_{kkii}-\phi_{kiik})+\phi P_{kl}(\phi_{ilki}-\phi_{klii})+\phi P_{kl}(\phi_{kiil}-\phi_{iikl})\notag\\
=&\phi P_{nn}(-\phi_{mk}R_{mk}-\phi_m R_{mk,k})+\phi P_{kl}(\phi_{mi}R_{mlik}+\phi_m
R_{ikml,i})+\phi P_{kl}{(\phi_{ml}R_{mk}+\phi_mR_{mk,l})}\notag
\end{align}
where we have used the Ricci identity in the last equality.

Substituting the following identities into \eqref{eqn:mid2} and \eqref{eqn:mid3},
\begin{equation*}
R_{ij}=(n-2)P_{ij}+\frac{R}{2(n-1)}g_{ij},\quad R_{ij,i}=\frac{R_j}{2},\quad
P_{kk}=\frac{R}{2(n-1)};
\end{equation*}
\begin{equation*}
R_{ijkl}=W_{ijkl}+P_{ik}g_{jl}+P_{jl}g_{ik}-P_{il}g_{jk}-P_{jk}g_{ik}, \;
P_{kk,i}=P_{ik,k};
\end{equation*}
\begin{equation*}\begin{split}
R_{ikml,i}&=R_{kl,m}-R_{km,l}
=(n-2)C_{klm}+\frac{\nabla_m R g_{kl}-\nabla_l R g_{km}}{2(n-1)}.
\end{split}\end{equation*}
after making some cancelations we see that the left hand side of \eqref{eqn:mid} becomes
\begin{align}\label{eqn:1}
&\int \phi\delta^{mnj}_{kli}P_{km}\phi_{nl}\phi_{ij}\\
\notag =&\int \Big[-\phi\phi_{kl}C_{lki,i}+\frac{4-n}{2(n-1)}\phi R\phi_{mk}P_{mk}-\frac{R^2}{4(n-1)^2}\phi\Delta\phi-\frac{n-2}{4(n-1)^2}\phi R \phi_m R_{,m}\\
  &\quad +\phi P_{kl}\phi_{mi}W_{mlik}+\phi |P_{kl}|^2\Delta \phi+(n-4)\phi P_{kl}\phi_{ml}P_{mk}+n\phi P_{kl}\phi_m P_{kl,m}\notag\\
  &\quad +2\phi P_{km,i}\phi_j W_{jmik}-2\phi P_{km,i}\phi_k P_{mi}\Big].\notag
\end{align}
On the other hand, by divergence theorem, we see that the other three terms on the last
of \eqref{eqn:2nd} are
\begin{align}\label{eqn:2}
\frac{2}{3}\int P_{ij}\phi_k\phi_l W_{ikjl}=&
   \frac{2}{3}\int -\phi P_{ij,k}\phi_l W_{ikjl}-\phi P_{ij}\phi_{kl}W_{ikjl}-\phi P_{ij}\phi_l W_{ikjl,k}\\
   =&\int -\frac{2}{3}\phi P_{ij,k}\phi_l W_{ikjl}-\frac{2}{3}\phi P_{ij}\phi_{kl}W_{ikjl}
   -\frac{2(n-3)}{3}\phi P_{ij}\phi_l C_{ijl},\notag
\end{align}

\begin{equation}\label{eqn:3}\begin{split}
-\frac{4}{3}\int \phi_k\phi_{ij}C_{ijk}
                =\int
                \frac{4}{3}\phi\phi_{ijk}C_{ijk}+\frac{4}{3}\phi\phi_{ij}C_{ijk,k},
\end{split}\end{equation}

\begin{align}\label{eqn:4}
&-\int (n-2)|\nabla\phi|^2\sigma_2(g)\\ =&\int (n-2)\phi\Delta\phi\sigma_2(g)
+(n-2)\phi \phi_i(\sigma_2(g))_{,i}\notag\\
=&\int \frac{(n-2)}{2}\phi\Delta\phi \Big(\frac{R^2}{4(n-1)^2}-|P_{kl}|^2\Big)
+(n-2)\phi\phi_{i}\Big(\frac{RR_{,i}}{4(n-1)^2}-P_{kl}P_{kl,i}\Big)\notag\\
=&\int \frac{(n-2)\phi\Delta\phi R^2}{8(n-1)^2}-\frac{n-2}{2}\phi\Delta\phi |P_{kl}|^2 +
\frac{(n-2)\phi R\phi_i R_{,i}}{4{(n-1)}^2}
                         -(n-2)\phi\phi_i P_{kl}P_{kl,i}.\notag
\end{align}
By combining equations \eqref{eqn:1}, \eqref{eqn:2}, \eqref{eqn:3} and \eqref{eqn:4} and
doing some cancelations,
 we conclude that the last four terms on the right hand side of \eqref{eqn:2nd} are equal to

\begin{align}\label{eqn:a}
-\frac{1}{8}&\int\Big[\frac{1}{3}\phi\phi_{kl}B_{kl}+(n-4)\phi
T_{2ij}\phi_{ij}-\frac{2(n-6)}{3}\phi
P_{ij}\phi_k C_{ijk}\\
  &\qquad\qquad  +\frac{4}{3}\phi C_{kmi}\phi_j W_{jmik}+\frac{4}{3}\phi \phi_{ijk}C_{ijk}\Big]dv,\notag
\end{align}
where we have used $T_{2ij}=\sigma_2(g)\delta_{ij}-\sigma_1(g)P_{ij}+P_{ik}P_{kj}$ and
$B_{ij}=C_{ijk,k}+P_{kl}W_{ikjl}$. Moreover,
\begin{equation*}\begin{split}
&\frac{4}{3}\phi C_{ijk}\phi_{jik}\\
=&\frac{4}{3}\phi C_{ijk}(\phi_{jki}+\phi_m R_{mjik})\\
=&\frac{4}{3}\phi C_{ijk}\phi_m(W_{mjik}+P_{mi}g_{jk}+P_{jk}g_{mi}-P_{mk}g_{ij}-P_{ij}g_{mk})\\
=&\frac{4}{3}\phi C_{ijk}\phi_m W_{mjik}+\frac{4}{3}\phi C_{ijj}\phi_m P_{mi}+\frac{4}{3}\phi C_{ijk}\phi_i P_{jk}-\frac{4}{3}\phi C_{iik}\phi_m P_{mk}-\frac{4}{3}\phi C_{ijk}P_{ij}\phi_k\\
=&-\frac{4}{3}\phi \phi_m C_{ijk}W_{mjki}-\frac{4}{3}\phi C_{ijk}P_{ij}\phi_k,
\end{split}
\end{equation*}
where we used $\sum_i C_{iik}=0$ and $C_{ijk}=-C_{ikj}$.

Thus it follows that \eqref{eqn:a} is equal to
\begin{equation}\label{eqn:mf}\begin{split}
-\frac{1}{8}\int\Big[\frac{1}{3}\phi\phi_{kl}B_{kl}+(n-4)\phi
\phi_{ij}T_{2ij}-\frac{2(n-6)}{3}\phi\phi_k P_{ij} C_{ijk}
    -\frac{4}{3}\phi\phi_k C_{ijk}P_{ij}\Big]dv.
\end{split}\end{equation}
Substituting \eqref{eqn:mf} into \eqref{eqn:2nd}, we conclude that
\begin{align}\label{eqn:2nd1}
&\frac{d^2}{dt^2}\Big|_{t=0} F(g_t) \\ =&\int \Bigg\{(n-6)\psi
v^{(6)}+(n-6)^2\phi^2v^{(6)}(g)-\frac{1}{8}\bigg[-2(n-5)\phi\phi_{ij} T_{2ij}
-\frac{2(n-5)}{3(n-4)}\phi \phi_{ij}B_{ij}\notag \\ & +\frac{4(n-5)}{3}\phi\phi_k
P_{ij}C_{ijk}+ \frac{1}{3}\phi\phi_{kl}B_{kl}+ (n-4)\phi
\phi_{ij}T_{2ij}-\frac{2(n-6)}{3}\phi\phi_k P_{ij} C_{ijk}
    -\frac{4}{3}\phi\phi_k C_{ijk}P_{ij}\bigg]\Bigg\}\notag \\
=&\int\Bigg\{ (n-6)\psi v^{(6)}+(n-6)^2\phi^2v^{(6)}+\frac{1}{8}\bigg[(n-6)\phi
\phi_{ij}T_{2ij}+\frac{n-6}{3(n-4)}\phi\phi_{ij} B_{ij}\notag\\ & \qquad
    - \frac{2(n-6)}{3}\phi\phi_k P_{ij} C_{ijk}\bigg]\Bigg\}\notag.
\end{align}

\noindent{\bf Proof of Theorem 1.2}  By Theorem \ref{thm:thm2}, at the critical metric of
the functional $\mathcal{F}_3(g)$, it holds that $\vs(g)$ should be constant, and it
follows that $F(g)= V\vs(g)$. By our notations
 $\mathcal{F}_3[g_t]=\frac{F(g_t)}{(\int\dv)^{(n-6)/n}}$. By use of \eqref{eqn:2nd1}
 and \eqref{eqn:main1}, at
 the critical metric $g$, we have 
\begin{align}\label{eqn:final}
&\frac{d^2}{dt^2}\Big|_{t=0} \mathcal{F}_3[g_t]\\
=& \frac{d^2}{dt^2}\Big|_{t=0} F(g_t)\cdot V^{-\frac{n-6}{n}} + 2 \frac{d}{dt}\Big|_{t=0}
F(g_t)\cdot
\frac{d}{dt}\Big|_{t=0}V(g_t)^{-\frac{n-6}{n}}+F(g)\frac{d^2}{dt^2}\Big|_{t=0}V(g_t)
^{-\frac{n-6}{n}}\notag\\
=& \frac{d^2}{dt^2}\Big|_{t=0} F(g_t)\cdot V^{-\frac{n-6}{n}}
    -\frac{2(n-6)}{n}\bigg(\frac{d}{dt}\Big|_{t=0}F(g_t)\bigg)\cdot V^{-\frac{(2n-6)}{n}}\int n\phi \notag\\ &\quad + F\Bigg\{
               (n-6)(2n-6)V^{-\frac{(3n-6)}{n}}\Big(\int \phi\Big)^2 -n(n-6)V^{-\frac{(2n-6)}{n}}\int \phi^2-
                   (n-6)V^{-\frac{(2n-6)}{n}}\int \psi\Bigg\}\notag\\
=&V^{-\frac{n-6}{n}}\Bigg\{\frac{d^2}{dt^2}\Big|_{t=0}
F(g_t)-2(n-6)^2v^{(6)}(g)V^{-1}\Big(\int \phi\Big)^2
    +(n-6)(2n-6)v^{(6)}(g)V^{-1}\Big(\int \phi\Big)^2\notag\\ &\qquad -n(n-6)v^{(6)}(g)\int \phi^2-
        (n-6)v^{(6)}(g)\int \psi \Bigg\}\notag\\
=&V^{-\frac{n-6}{n}}\Bigg\{\frac{d^2}{dt^2}\Big|_{t=0} F(g_t) +
6(n-6)v^{(6)}(g)V^{-1}\Big(\int \phi\Big)^2 - n(n-6)v^{(6)}(g)\int \phi^2
                                 -(n-6)v^{(6)}(g)\int \psi\Bigg\}\notag\\
=&V^{-\frac{n-6}{n}}\Bigg\{\int \bigg[-6\vs(g)\Big(\phi-V^{-1}\int \phi\Big)^2
+\frac{\phi\phi_{kl}}{24(n-4)}B_{kl}+\frac{1}{8}
 \phi\phi_{mk}T_{2mk}\notag\\
  &\qquad\qquad - \frac{1}{12}\phi C_{ijk}P_{ij}\phi_k\bigg]dv \Bigg\}\notag.
  \end{align} If we define an operator $\mathcal L$ by
\begin{equation*}
\mathcal{L}(f):=\frac{B_{ij}f_{ij}}{24(n-4)}+\frac{1}{8}T_{2ij}f_{ij}-
\frac{1}{12}P_{ij}C_{ijk}f_k,
\end{equation*}
for $f\in C^\infty(M)$. It is easy to see that $\mathcal{L}$ is self-adjoint with respect
to the $L^2$ inner product of the Riemannian manifold. Indeed,  for any two smooth
functions $f$ and $h$, we have
\begin{equation*}\begin{split}
\langle\mathcal{L}(f),h\big\rangle&=\int_M \mathcal{L}(f)h dv\\
&=\int_M \Big[\frac{B_{ij}f_{ij}h}{24(n-4)}+\frac{1}{8}T_{2ij}hf_{ij}-
\frac{1}{12}P_{ij}C_{ijk}f_kh\Big]dv\\
&=\int_M \Big[-\frac{B_{ij}f_ih_j}{24(n-4)}-\frac{1}{8}T_{2ij}f_i h_j-\frac{B_{ij,j}f_i
h}{24(n-4)}-\frac{1}{8}T_{2ij,j}f_i h-\frac{1}{12}P_{ij}C_{ijk}f_k h\Big]dv\\
 &=\int_M \Big[-\frac{B_{ij}f_ih_j}{24(n-4)}-\frac{T_{2ij}f_ih_j}{8}\Big] dv\\
&=\langle f,\mathcal{L}(h)\rangle, \end{split}\end{equation*} where we have used (2)and (3) in Lemma
\ref{lem:lem1}, \eqref{eqn:tij} and integration by parts. Denote $\phi-V^{-1}\int \phi$ by $\bar{\phi}$. From
\eqref{eqn:final}, we see that
\begin{align}
&\frac{d^2}{dt^2}\Big|_{t=0} \mathcal{F}_3(g_t)\\
&=(n-6)V^{-\frac{n-6}{n}}\Bigg\{\int \Big[-6\vs(g)\bar{\phi}^2+\mathcal{L}(\phi)\phi   \Big]dv\Bigg\}\notag\\
&=(n-6)V^{-\frac{n-6}{n}}\Bigg\{\int
\Big[-6\vs(g)\bar{\phi}^2+\mathcal{L}(\bar{\phi})\bar{\phi} \Big]dv\Bigg\}.\notag
\end{align}
Thus we complete the proof of Theorem \ref{thm:thm1}.\qed

 To prove Theorem \ref{thm:main}, we need the
following famous theorem.
\begin{thm}[Lichnerowicz and Obata, see e.g. \cite{LI}]\label{thm:thm3}
Let $M$ be an $n$-dimensional compact manifold. Suppose the Ricci curvature of $M$ is
bounded from below by
\begin{equation*}Ric\ge (n-1)K\end{equation*} for some positive constant $K$, then the first nonzero
eigenvalue of the Laplacian on $M$ must satisfy  $$\lambda_1\ge nK.$$ Moreover, equality
holds if and only if $M$ is isometric to a standard sphere of radius
$\frac{1}{\sqrt{K}}$.
\end{thm}
By the min-max principle, for the first nonzero eigenvalue $\lambda_1$ of Laplacian, it
holds that
\begin{equation}\label{eqn:cha}
\lambda_1\int_M f^2dv\le \int_M|\nabla f|^2dv,
\end{equation}
for any $f\in C^\infty(M)$ satisfying $\int_M f dv=0$.

\noindent\medskip {\bf Proof of Theorem 1.3}. Note that an Einstein manifold $(M^n,g)$ is
a critical metric in $[g]$, i.e. it satifies (1.2).
 Now let $(M^n,g)$ be an Einstein manifold with positive scalar
curvature, then it follows from Theorem \ref{thm:thm3} and \eqref{eqn:cha} that
\begin{equation}\label{eqn:res}
\frac{R}{n-1}\int_M \bar{\phi}^2dv\le \int_M|\nabla \bar{\phi}|^2dv.
\end{equation}
Note that for an Einstein manifold, $\vs(g)=-\frac{(n-2)R^3}{386n^2(n-1)^2}$,
$\mathcal{L}(\phi)= \frac{(n-2)R^2}{64n^2(n-1)}\Delta \phi$. Hence, we see that
\begin{align}\label{eqn:424}
\frac{d^2}{dt^2}\Big|_{t=0}\mathcal{F}_3[g_t]&=(n-6)V^{-\frac{n-6}{n}}\int_M
\Big[\frac{(n-2)R^3}{64n^2(n-1)^2}\bar{\phi}^2 + \frac{(n-2)R^2\bar{\phi}}{64n^2(n-1)}\Delta\bar{\phi}\Big]dv ,\notag\\
&=\frac{(n-2)(n-6)R^2}{64n^2(n-1)}V^{-\frac{n-6}{n}}\int_M\Big[\frac{R\bar{\phi}^2}{n-1}
      -|\nabla\bar{\phi}|^2\Big]dv\notag\\
&\le\frac{(n-2)(n-6)R^2}{64n^2(n-1)}V^{-\frac{n-6}{n}} \Big(\frac{R}{n-1}-\lambda_1\Big)\int_M \bar{\phi}^2 dv\\
&\le 0,\notag
\end{align}
with equality holds if and only if $\lambda_1=\frac{R}{n-1}$. Hence, by Theorem
\ref{thm:thm3}, in this case $(M,g)$ is isometric to the standard sphere $S^n$.

Therefore, we prove that an Einstein manifold with positive scalar curvature must be a
strict local maximum ``point'' within its conformal class $[g]$ unless $(M,g)$ is
isometric to $S^n$ with a multiple of the standard metric. We complete the proof of
theorem 1.3. \qed

\noindent\medskip

\begin{rem}
Let $(M^n,g)$ be an $n$-dimensional Einstein manifold with nonpositive scalar curvature,
then we have from the proof of Theorem \ref{thm:main} (see \eqref{eqn:424})
$$
\frac{d^2}{dt^2}\Big|_{t=0}\mathcal{F}_3(g_t)\le 0,
$$
that is, it is stable.
\end{rem}
\begin{rem}
When $M^n$ is an Einstein manifold with positive scalar curvature with dimension $n=5$,
we see from \eqref{eqn:424} that
\begin{equation}\label{eqn:test}
\frac{d^2}{dt^2}\Big|_{t=0}\mathcal{F}_3[g_t]\ge 0,
\end{equation}
with equality if and only if $\lambda_1=\frac{R}{4}$. Theorem \ref{thm:thm3} shows that
in this case $(M^5,g)$ is isometric to the sphere $S^5$ with the standard metric up to a
multiple of constant. And we see that this Einstein metric is a strict local minimum of
the functional $\mathcal{F}_3$ within its conformal class if the equality does not hold
in \eqref{eqn:test}.
\end{rem}

\begin{rem}
Let $\mathcal{T}_{ij}(g)=T_{2ij}(g)+\frac{1}{n-4}(B_g)_{ij}$, we have
\begin{equation*}
\sum_j\nabla^{j}\mathcal{T}_{ij}=0,
\end{equation*}
that is, $\mathcal{T}_{ij}$ is a divergence-free tensor. We observe that
$v^{(6)}(g)=-24\sum_{ij}\mathcal{T}_{ij}(g)(P_g)_{ij}$.
\end{rem}

\noindent Bin Guo: {\sc Department of Mathematical Sciences, Tsinghua University, Beijing 100084, People's
Republic of China}\ \  Email: guob07@mails.tsinghua.edu.cn

\noindent Haizhong Li: {\sc Department of Mathematical Sciences, Tsinghua University, Beijing 100084,
People's Republic of China} \ \ E-mail: hli@math.tsinghua.edu.cn

\begin{thebibliography}{99}


\bibitem[CF]{CF08}
S. -Y. A. Chang and H. Fang, {A class of variational functionals in conformal geometry},
Int. Math. Res. Not. (2008), No. 7, rnn008, 16 pages, arXiv: math/0610773.


\bibitem[G]{G00}
C. Robin Graham, {Extended obstruction tensors and renormalized volume coefficients},
Advances in Math. , 220(2009), no.6, 1956-1985.

\bibitem[GHL]{GHL}Bin Guo, Zheng-Chao Han and Haizhong Li, Two Kazdan-Warner type identities for the renormalized volume coefficients and
Gauss-Bonnet curvatures of a Riemannian metric, arXiv: math/0911.4649.
\bibitem[GJ]{GJ07}
C. Robin Graham and A. Juhl, {Holographic formula for $Q$ curvature}, Advances in Math. ,
216(2007), 841-53.

\bibitem[L]{LI} Peter Li, Lecture notes on geometric analysis, Lecture Notes Series No. 6 - Research Institute
 of Mathematics and
Global Analysis Research Center, Seoul National University, Seoul, 1993.

\bibitem[R]{Re} R. Reilly, { Applications of the Hessian operator in
a  Riemannian manifold}, Indiana Univ. Math. J., 26(3)(1977), 459-472.

\bibitem[V]{V} J. Viaclovsky, Conformal geometry, contact geometry, and the calculus of variations,
Duke Math. J. 101, No. 2 (2000), 283-316.
\end{thebibliography}
\end{document}